\documentclass{article}
\usepackage{style}

\title{On the Lie enveloping algebra of a post-Lie algebra}

\author{
Kurusch Ebrahimi-Fard\footnote{\small ICMAT, Calle Nicol\'as Cabrera 13-15, Campus de Cantoblanco, UAM, Madrid, Spain. On leave from Université de Haute Alsace, 18 Rue des Frères Lumière, 68093 Mulhouse Cedex, France. {\small{kurusch@icmat.es}}} 
\and 
Alexander Lundervold\footnote{\small Dept. of Computing, Mathematics and Physics, Bergen University College, Postboks 7030, 5020 Bergen, Norway. {\small{alexander.lundervold@hib.no}}}
\and 
Hans Munthe-Kaas\footnote{\small Dept. of Mathematics, University of Bergen, Postbox 7800
N-5020 Bergen, Norway. {\small{hans.munthe-kaas@math.uib.no}}}}

\begin{document}

\maketitle

\begin{abstract}
We consider pairs of Lie algebras $\g$ and $\bar{\g}$, defined over a common vector space, where the Lie brackets of $\g$ and $\bar{\g}$ are related via a post-Lie algebra structure. The latter can be extended to the Lie enveloping algebra $\mathcal{U}(\g)$. This permits us to define another associative product on $\mathcal{U}(\g)$, which gives rise to a Hopf algebra isomorphism between $\mathcal{U}(\bar{\g})$ and a new Hopf algebra assembled from $\mathcal{U}(\g)$ with the new product.

For the free post-Lie algebra these constructions provide a refined understanding of a fundamental Hopf algebra appearing in the theory of numerical integration methods for differential equations on manifolds. In the pre-Lie setting, the algebraic point of view developed here also provides a concise way to develop Butcher's order theory for Runge--Kutta methods.    
\end{abstract}

.


\section{Introduction}
\label{sect:intro}

Classical numerical integrators aim at approximating flows of vector fields given by differential equations 
\[
	\dot{y}(t) = F(y(t)), \quad y(0) = y_0,
\]
where $F$ is a vector field on $\RR^n$. The so-called \emph{Lie-group integrators} are generalizations to differential equations evolving on manifolds $\mathcal{M}$. Given a Lie group $G$ acting in a transitive manner on $\mathcal{M}$, the Lie group integrators approximate differential equations written in the form
\begin{equation}\label{eq:lgdiffeqn}
	\dot{y}(t) = F(y(t)) =  f(y(t))\cdot y(t), \quad y(0) = y_0,
\end{equation}	
where $f \colon \mathcal{M} \rightarrow \g$ represents the vector field $F$ via a map to $\g$, the Lie algebra of $G$. The basic problem of numerical Lie group integration is the approximation of the exact solution by computations in $\g$, the exponential $\exp\colon \g\rightarrow G$ and the action of $G$ on $\mathcal{M}$. See \cite{iserles2000lgm} for a comprehensive review.

Answering questions related to order conditions for Lie group integrators relies on an understanding of the algebraic structure of non-commuting vector fields on $\mathcal{M}$ generated from the vector field $f$. Our work addresses such algebraic and combinatorial aspects underlying the theory of Lie--Butcher series \cite{munthe-kaas1995lbt}. More precisely, we describe the setting of an important commutative Hopf algebra of planar rooted trees \cite{munthe-kaas2008oth}, and the corresponding group and Lie algebra of characters and infinitesimal characters, respectively. In a nutshell, the theory of Lie--Butcher series results from merging Lie theory with the theory of Butcher's $B$-series. The latter is well-known, and plays a prominent role in the theory of numerical integration for differential equations on vector spaces \cite{butcher1972aat,butcher2008nmf,hairer2006gni}. 

The underlying context is the following: We have a Lie algebra $\g=(\mathcal{V}, [-, -])$ over a field $k$ of characteristic zero. In addition, the vector space $\mathcal{V}$ is equipped with a binary product $\tr: \mathcal{V} \otimes \mathcal{V} \to \mathcal{V}$, such that the bracket 
\begin{equation*}
	\llbracket x,y \rrbracket := x \tr y - y \tr x + [x,y]
\end{equation*}
defines a second Lie algebra $\bar{\g}=(\mathcal{V}, \llbracket -,- \rrbracket )$ on $\mathcal{V}$. The Lie algebra $\g$ together with the appropriate relations between its Lie bracket and the product $\tr: \mathcal{V} \otimes \mathcal{V} \to \mathcal{V}$ define what is called a post-Lie algebra. It reduces to a pre-Lie algebra (also called a \emph{Vinberg algebra}) \cite{cartier2010val,manchon2011ass} if the Lie algebra $\g$ is abelian. 

Post-Lie algebras appear in both geometry and algebra \cite{munthe-kaas2013opl,vallette2007hog}. From a geometric viewpoint, they can be introduced by means of an archetypal example from differential geometry\footnote{To keep the discussion as simple as possible, we discuss an affine connection on the tangent bundle. A similar construction is also valid on more general vector bundles such as the Atiyah Lie algebroid~\cite{munthe-kaas2013opl}.}. 

Recall that an {affine connection} (covariant derivative) on the space of smooth vector fields $\mathcal{X}(\mathcal{M})$ on a manifold $\mathcal{M}$ is a map $\nabla \colon \mathcal{X}(\mathcal{M})\times \mathcal{X}(\mathcal{M}) \rightarrow\mathcal{X}(\mathcal{M})$ satisfying $\nabla_{fx} y = f \nabla_x y$ and $\nabla_x (fy) = df(x) y + f\nabla_x y$ for any $x,y\in \mathcal{X}(\mathcal{M})$ and scalar field $f: \mathcal{M}\rightarrow \mathbb{R}$. It yields a (non-commutative and non-associative) $\mathbb{R}$-bilinear product on $\mathcal{X}(\mathcal{M})$, which we denote $x \tr  y := \nabla_x y.$ The torsion of the connection is a skew-symmetric tensor $\mathrm{T} \colon T\mathcal{M} \wedge T\mathcal{M} \rightarrow T\mathcal{M}$ defined as
\begin{equation}
	\mathrm{T}(x,y) = x\tr y -y\tr x - \llbracket x,y\rrbracket_J ,
	\label{eq:torsion}
\end{equation}
where $\llbracket\cdot,\cdot\rrbracket_J$ denotes the Jacobi--Lie bracket of vector fields, defined such that  $\llbracket x,y\rrbracket_J(\phi) = x(y(\phi))-y(x(\phi))$ for all vector fields $x,y$ and scalar fields $\phi$. The curvature tensor $\mathrm{R}\colon T\mathcal{M} \wedge T\mathcal{M} \rightarrow \mbox{End}(T\mathcal{M})$ is defined as
\begin{align}
	\mathrm{R}(x,y)z &= x\tr (y\tr z)- y\tr (x\tr z)- {\llbracket x,y\rrbracket_J}\tr z  \\
				  &= a_\tr(x,y,z)-a_\tr(y,x,z)+\mathrm{T}(x,y)\tr z,
	\label{eq:curvature}
\end{align}
where $a_\tr(x,y,z):= x\tr (y\tr z)-(x\tr y)\tr z$ is the associator with respect to the product $\tr$. The relationship between torsion and curvature is given by the Bianchi identities
\begin{align}
	\sum_{\circlearrowleft}(\mathrm{T}(\mathrm{T}(x,y),z)+(\nabla_x \mathrm{T})(y,z)) 
	&= \sum_{\circlearrowleft}(\mathrm{R}(x,y)z)\label{eq:bianchi1}\\
	\sum_{\circlearrowleft}((\nabla_x\mathrm{R})(y,z)+\mathrm{R}(\mathrm{T}(x,y),z)) 
	&= 0 ,\label{eq:bianchi2}
\end{align}
where $\sum_{\circlearrowleft}$ denotes the sum over the three cyclic permutations of $(x,y,z)$. If a connection is flat $\mathrm{R}=0$ and has constant torsion $\nabla_x \mathrm{T}=0$, then \eqref{eq:bianchi1} reduces to the Jacobi identity $\sum_{\circlearrowleft}(\mathrm{T}(\mathrm{T}(x,y),z)) = 0$. Hence the torsion defines a Lie bracket $[x,y] := -\mathrm{T}(x,y)$, which is related to the Jacobi--Lie bracket by \eqref{eq:torsion}. The covariant derivation formula $\nabla_x (\mathrm{T}(y,z)) = (\nabla_x \mathrm{T})(y,z)+\mathrm{T}(\nabla_x y,z)+\mathrm{T}(y,\nabla_x z)$ together with $\nabla_x \mathrm{T}=0$ then imply
\begin{equation}
	x\tr [y,z] = [x\tr y, z] + [y,x\tr z] .\label{eq:postlie1}
\end{equation}
On the other hand, \eqref{eq:curvature} together with $\mathrm{R}=0$ imply
\begin{equation}
	[x,y]\tr z = a_\tr(x,y,z) - a_\tr(y,x,z) .\label{eq:postlie2}
\end{equation}
In Section \ref{sect:PrePostLie} relations \eqref{eq:postlie1} and \eqref{eq:postlie2} are formalized to the notion of post-Lie algebras.

Note that for a connection which is both flat, $\mathrm{R}=0$, and torsion free, $\mathrm{T}=0$, equation \eqref{eq:curvature} implies that $a_\tr(x,y,z) - a_\tr(y,x,z)=0$, and we obtain a (left) pre-Lie algebra. 
%
The free pre-Lie algebra can be described as the space of non-planar rooted trees with product given by grafting of trees~\cite{chapoton2001pla}. We remark that already in the 1850s Arthur Cayley \cite{cayley1857ott} considered rooted trees as a representation of combinatorial structures related to the free pre-Lie algebra. More than a century later, these structures form the foundation of John Butcher's theory of $B$-series~\cite{butcher1963cft,butcher1972aat}, which has become an indispensable tool in the analysis of numerical integration~\cite{hairer2006gni}.

Pre-Lie structures on non-planar rooted trees lead to Hopf algebras of combinatorial nature. Combinatorial Lie and Hopf algebras, in particular those defined on rooted trees, have recently attracted a great deal of attention \cite{connes1998har, grossman1989has, hoffman2003cor, murua2006tha}. 

The basic setting for Lie--Butcher series is provided by a combinatorial Hopf algebra on planar rooted trees, which accompanies a post-Lie structure on the trees \cite{munthe-kaas2008oth}. One of the main goals of this work is to provide a precise description of this connection in the context of the free post-Lie algebra. Our approach follows Guin and Oudom \cite{oudom2008otl} by extending the post-Lie structure on a Lie algebra $\g$ to the corresponding Lie enveloping algebra $\mathcal{U}(\g)$. This permits us to define another associative product on $\mathcal{U}(\g)$, compatible with its usual coalgebra structure. The Hopf algebra assembled from $\mathcal{U}(\g)$ and the new product is isomorphic to the Hopf algebra $\mathcal{U}(\bar{\g})$, where $\bar{\g}$ is the Lie algebra defined over the same vector space as $\g$, whose Lie bracket is defined in terms of a post-Lie algebra structure. 

Once the combinatorial Hopf algebra for Lie--Butcher series has been unfolded we explore some applications, in particular to the order theory of numerical methods on manifolds.

\medskip
 
The outline of the paper is as follows. In the next section we recall the definition of pre- and post-Lie algebra. Section \ref{sect:UnivLieEnvAlgPostLie} generalizes some important result from \cite{oudom2008otl}, i.e.~the extension of post-Lie structures to universal enveloping Lie algebras, and a Hopf algebra isomorphism between two enveloping algebras. In Section \ref{sect:RK} we study numerical integration on pre- and post-Lie algebras.


\section{Pre- and post-Lie algebras}
\label{sect:PrePostLie}

We begin by defining pre-Lie algebras, see \cite{burde2006left,cartier2010val,manchon2011ass} for further details. A field $k$ of characteristic zero,  e.g.~$k \in \{\mathbb{R},\mathbb{C}\}$, is fixed once and for all. Let $\mathcal{P}$ be a vector space equipped with a bilinear product $\tr : \mathcal{P} \times \mathcal{P} \to \mathcal{P}$, satisfying the left pre-Lie relation
\begin{equation}
\label{def:pL}
(x\tr y)\tr z - x\tr (y\tr z) = (y\tr x)\tr z - y\tr (x\tr z).
\end{equation}
We call $(\mathcal{P},\tr)$ a left pre-Lie algebra. Note that identity \eqref{def:pL} can be written as
\begin{equation}
	\ell_{[x,y]\tr}(z)=[\ell_{x \tr},\ell_{y\tr}](z),
\end{equation}
where the linear map $\ell_{x \tr}: \mathcal{P} \to \mathcal{P}$ is defined by $\ell_{x\tr}y:=x\tr y$. The bracket on the left-hand side is defined by $[x,y]:=x \tr y - y \tr x$. As a consequence this bracket satisfies the Jacobi identity, turning $\mathcal{P}$ into a Lie algebra.  

The notion of post-Lie algebras was discovered from two distinct directions \cite{munthe-kaas2008oth,vallette2007hog}. We follow the definition given in \cite{munthe-kaas2013opl}. The \emph{associator} is defined by $a_{\tr}(x,y,z) :=  x \tr (y \tr z) - (x \tr y) \tr z.$

\begin{Definition}  {\rm{\cite{munthe-kaas2013opl}}} \label{def:post-Lie}
A post-Lie algebra $(\mathcal{A}, [-, -], \tr)$ is a Lie algebra $\g=(\mathcal{A}, [-, -])$ together with a bilinear product $\tr: \mathcal{A} \times \mathcal{A} \to \mathcal{A}$ such that for all $x,y,z \in \mathcal{A}$
\begin{align}
	x \tr [y,z] &= [x \tr y,z] + [y,x \tr z] \label{def:post-lie1}\\
	{[x,y]} \tr z &= a_{\tr}(x,y,z) - a_{\tr}(y,x,z).\label{def:post-lie2}
\end{align}
\end{Definition}

\begin{Proposition} \label{prop:post-lie}
Let $(\mathcal{A}, [-, -], \tr)$ be a post-Lie algebra. The bracket
\begin{equation}
\label{def:post-lie3}
	\llbracket x,y \rrbracket = x \tr y - y \tr x + [x,y]
\end{equation}
satisfies the Jacobi identity for all $x, y \in \mathcal{A}$. 
\end{Proposition}

\begin{Remark} 
If $(\mathcal{A}, [-, -])$ is an abelian Lie algebra, then $(\mathcal{A}, \tr)$ reduces to a left pre-Lie algebra with corresponding Lie bracket \eqref{def:post-lie3}.
\end{Remark}

The post-Lie algebra defines relations between two Lie algebras $[-,-]$ and $\llbracket-,-\rrbracket$, which will be explored in the sequel. It should be noted that the role of these two brackets is not at all symmetric in the definition of the post-Lie structure. Whereas the above definition relating $[-,-]$ and $\tr$ seems, algebraically speaking, to be the nicest definition, we
could instead follow a more geometrically motivated definition, where we equip a general Lie algebra $\{{\mathcal A},\llbracket -,-\rrbracket\}$  with a product $\tr\colon {\mathcal A}\times{\mathcal A}\rightarrow {\mathcal A}$ such that the `curvature' of $\tr$ is zero:
\begin{equation}\label{eq:R0}R(x,y)z := x\tr(y\tr z) - y\tr(x\tr z) - \llbracket x,y\rrbracket\tr z = 0 .\end{equation}
We then define the `torsion' of $\tr$ as
\begin{equation}[x,y] := \llbracket x,y \rrbracket + y\tr x -x\tr y, \end{equation}
and require that this satisfies the condition 
\begin{equation}\label{eq:nT0}
x\tr [y,z] = [x\tr y,z] + [y,x\tr z]
\end{equation}
(geometrically, $\nabla T=0$).

\begin{Proposition} Let $\{{\mathcal A},\llbracket\cdot,\cdot\rrbracket\}$ be a Lie algebra with a product $\tr\colon {\mathcal A}\times{\mathcal A}\rightarrow {\mathcal A}$ 
satisfying~\eqref{eq:R0}--\eqref{eq:nT0}. Then $\{A,[\cdot,\cdot],\tr\}$ is post-Lie.
\end{Proposition}

\begin{proof}A direct computation shows that $[\cdot,\cdot]$ satisfies the Jacobi rule for a Lie bracket. The other conditions of post-Lie are baked into the defining equations.
\end{proof}

Post-Lie algebras always come in pairs called \emph{adjoint} post-Lie algebras.  E.g.\ in the case of vector fields on Lie groups, the algebras obtained from  left- and right trivialization of the vector fields are adjoint.  

\begin{Proposition}\label{prop:twist}
Let $(A,[\cdot,\cdot],\tr)$ be a post-Lie algebra and define the product $\vartr$ as
\begin{equation}
\label{eq:rightpostlie}
	x \vartr y := x \tr y + [x,y].
\end{equation}
Then the adjoint $(A,-[\cdot,\cdot],\vartr)$ of $(A,[\cdot,\cdot],\tr)$ is also a post-Lie algebra. The adjoint operation is an involution. The two adjoint post-Lie algebras generate the same `double bracket' $\llbracket \cdot,\cdot\rrbracket$.
\end{Proposition}


\section{Universal Lie enveloping algebras of a post-Lie algebra}
\label{sect:UnivLieEnvAlgPostLie}

Recall that we consider a post-Lie algebra as two Lie algebras $\g:=(\mathcal{A}, [\cdot, \cdot])$ and $\bar{\g}:=(\mathcal{A}, \llbracket \cdot, \cdot \rrbracket)$ defined over the vector space $\mathcal{A}$, related via the post-Lie product \eqref{def:post-lie3} in Proposition \ref{prop:post-lie}
\[
	\llbracket x,y \rrbracket = x \tr y - y \tr x + [x,y].
\]

It is natural to explore this relation on the level of universal Lie enveloping algebras. This is analogous to the approach in \cite{oudom2008otl}, where the universal Lie enveloping algebra of a pre-Lie algebra was studied. In \cite{oudom2008otl} it was shown that the universal Lie enveloping algebra of the Lie algebra obtained by antisymmetrization of a pre-Lie algebra $\mathcal{\lambda}$ is isomorphic as a Hopf algebra to the symmetric algebra $\mathcal{S}(\mathcal{\lambda})$ on $\mathcal{\lambda}$, equipped with a certain associative product $*: \mathcal{S}(\mathcal{\lambda}) \otimes \mathcal{S}(\mathcal{\lambda}) \to \mathcal{S}(\mathcal{\lambda})$ defined using the pre-Lie structure. 

Let $\mathcal{P}:=(\mathcal{A}, [\cdot, \cdot], \tr)$ be a post-Lie algebra, and $\mathcal{U}(\g)$ the universal enveloping algebra of the Lie algebra  $\g:=(\mathcal{A}, [\cdot, \cdot])$. We generalize the results of \cite{oudom2008otl} to universal Lie enveloping algebras built on Lie algebras related via the post-Lie relation \eqref{def:post-lie3}. 

The post-Lie product can be extended to a product on $\mathcal{U}(\g)$. We first extend it to a product mapping $\g \otimes \mathcal{U}(\g) \rightarrow \mathcal{U}(\g)$ and then, in Proposition \ref{prop:post-lieUg}, to all of $\mathcal{U}(\g)$. Let $x, t_1,\ldots, t_n \in \g$, and define
\begin{equation}
\label{def:trUg}
	x \tr \un := 0 \qquad x \tr (t_1\cdots t_n) := \sum_{i=1}^n t_1 \cdots t_{i-1}(x \tr t_i) t_{i+1} \cdots t_n.
\end{equation}
Note that this extension defines a natural action of $\mathcal P$ on $\mathcal{U}(\g)$, satisfying
\begin{align*}
	[x,y] \tr K &=    x \tr (y \tr K) - (x \tr y) \tr K - y \tr (x \tr K) + (y \tr x) \tr K\\
	x\tr [K_1,K_2] & = [x\tr K_1,K_2] + [K_1,x\tr K_2]
\end{align*}
for $K,K_1,K_2 \in \mathcal{U}(\g)$ and $x,y \in \g$, where $[K_1,K_2] = K_1K_2-K_2K_1$.

Recall that $\mathcal{U}(\g)$ with concatenation as product is a non-commutative, cocommutative Hopf algebra. The coshuffle coproduct is defined for $x \in \g \subset \mathcal{U}(\g)$ by $\Delta(x) := x \otimes \un + \un \otimes x$, and extended multiplicatively to all of $\mathcal{U}(\g)$. We use Sweedler's notation for the coproduct: $\Delta(T)=:T_{(1)} \otimes T_{(2)}$. The counit is denoted by $\epsilon: \mathcal{U}(\g) \to k$.

The remainder of this section contains statements and some proofs analogous to the ones in \cite[Section 2]{oudom2008otl}, slightly modified to match our setting. The omitted proofs can easily be deduced from \cite{oudom2008otl}.

\begin{Proposition} \label{prop:post-lieUg} 
Let $A,B,C \in \mathcal{U}(\g)$ and $x,y \in \g$. There is a unique extension of the post-Lie product $\tr$ from $\g$ to $\mathcal{U}(\g)$ given by
\begin{align}
	& \un \tr A  = A						\label{Ext1}\\
  	& xA \tr y 	= x \tr (A \tr y)	 - (x \tr A) \tr y 	\label{Ext2}\\
  	& A  \tr BC	= (A_{(1)} \tr B)(A_{(2)} \tr C).	        \label{Ext3}	 
\end{align}
\end{Proposition}
The result can be proven by induction as in \cite[Proposition 3.7]{oudom2008otl}: First, \eqref{Ext1} and \eqref{Ext3} together lead to the formula
\[
	x \tr (t_1 \cdots t_n) = \sum_{i=1}^n t_1 \cdots t_{i-1} (x\tr t_i )t_{i+1}  \cdots t_n,
\] 
and then \eqref{Ext3} and induction on the length of $B$ can be used to extend $A \tr B$ to monomials $A$ and $B$.

\begin{Proposition}\label{prop:relations}
Let $A,B,C \in \mathcal{U}(\g)$ and $x \in \g$. We have
\begin{align}
	& A \tr \un = \epsilon(A) \label{res1}\\
  	& \epsilon(A \tr B) = \epsilon(A)\epsilon(B) \label{res2}\\
  	& \Delta(A \tr B) = (A_{(1)} \tr B_{(1)}) \otimes (A_{(2)} \tr  B_{(2)}) \label{res3}\\
  	& xA \tr B = x \tr (A \tr B) - (x \tr A) \tr B \label{res4}\\
  	& A \tr (B \tr C) = (A_{(1)}(A_{(2)} \tr B)) \tr C \label{res5}
\end{align}
\end{Proposition}

\begin{proof}

The proofs of \eqref{res1}, \eqref{res2}, \eqref{res3}, \eqref{res4} are straightforward adaptions of the proofs in \cite[Proposition 3.9]{oudom2008otl}. We prove Equation \eqref{res5} by induction on the length of $A$. 
{\small
\begin{align}
	& xA \tr (B \tr C) = x \tr (A \tr (B \tr C)) - (x \tr A) \tr (B \tr C) \\
 	&\quad = x \tr(A_{(1)} (A_{(2)} \tr B) \tr C) - (x\tr A) \tr (B \tr C) \label{ind}\\
 	&\quad = (x A_{(1)} (A_{(2)} \tr B)) \tr C + (x \tr (A_{(1)}(A_{(2)} \tr B))) \tr C - (x\tr A)\tr (B\tr C)  \label{4}\\
 	&\quad = (x A_{(1)} (A_{(2)} \tr B)) \tr C + (A_{(1)} (x \tr (A_{(2)} \tr B))) \tr C \label{prev3}\\
 	& \qquad+ (x \tr A_{(1)})(A_{(2)} \tr B) \tr C - (x\tr A)\tr(B\tr C) \\
 	&\quad = (x A_{(1)} (A_{(2)} \tr B)) \tr C + (A_{(1)} (x (A_{(2)} \tr B))) \tr C \label{4-1}\\
 	& \qquad + ((A_{(1)} \tr A_{(2)}) \tr B) \tr C + (x\tr A_{(1)})(A_{(2)} \tr B) \tr C - (x \tr A)\tr(B\tr C)\\
	&\quad = (x A_{(1)} (A_{(2)} \tr B)) \tr C + (A_{(1)} (x A_{(2)} \tr B)) \tr C + \label{4-2}\\
 	& \qquad (A_{(1)} ((x \tr A_{(2)}) \tr B)) \tr C + (x \tr A_{(1)})(A_{(2)} \tr B) \tr C - (x\tr A)(B\tr C)
\end{align}
}
Here \eqref{ind} follows by induction, \eqref{4} from \eqref{res4}, \eqref{prev3} from \eqref{Ext3} in Proposition \ref{prop:post-lieUg}, \eqref{4-1} from \eqref{res4}, and \eqref{4-2} from \eqref{res4}.
We have
\begin{align}
	(x \tr A) \tr (B \tr C) 	&= (x \tr A)_{(1)} ((x \tr A)_{(2)} \tr B) \tr C \\
 					&= (x \tr A_{(1)})(A_{(2)} \tr B) \tr C + A_{(1)} ((x \tr A_{(2)}) \tr B) \tr C,
\end{align}
so that
\begin{align}
	xA \tr (B \tr C) 	&= (x A_{(1)} (A_{(2)} \tr B)) \tr C + (A_{(1)}(xA_{(2)} \tr B)) \tr C \\
 				&= (xA)_{(1)}((xA)_{(2)} \tr B)) \tr C.
\end{align}
\end{proof}

\begin{Proposition}\label{prop:postpre-lie}
On $(\mathcal{U}(\g),\tr)$ the product:
 \begin{equation}
\label{def:preLie}
	A * B := A_{(1)}(A_{(2)} \tr B)
\end{equation}
is associative. Moreover, $(\mathcal{U}(\g), *, \Delta)$ is a Hopf algebra.
\end{Proposition}

\begin{proof}
From \ref{prop:relations} we have
\begin{align}
	A * (B * C) &= A_{(1)} \left(A_{(2)} \tr ( B_{(1)} (B_{(2)} \tr C))\right) \\
 			&= A_{(1)} \left((A_{(2)} \tr B_{(1)}) (A_{(3)} \tr (B_{(2)} \tr C))\right) \label{ass1} \\
 			&= A_{(1)} \left((A_{(2)} \tr B_{(1)})((A_{(3)} (A_{(4)} \tr B_{(2)})) \tr C)\right) \label{ass2} \\
 			&= A_{(1)} \left((A_{(3)} \tr B_{(1)})((A_{(2)} (A_{(4)} \tr B_{(2)})))\tr C\right) \label{ass3} \\
 			&= (A_{(1)} (A_{(2)} \tr B)) * C \\
 			&= (A * B) * C,
\end{align}
where \eqref{ass1} follows from \eqref{Ext3}, \eqref{ass2} from \eqref{res5}, and \eqref{ass3} from cocommutativity. The compatibility between~\eqref{def:preLie} and the coproduct $\Delta$ follows from \eqref{Ext3} 
\end{proof}

Note that for $a$ and $b$ in $\g$, we have $\Delta(a)=a\otimes I + I\otimes a$, so
\[
	a * b = a \tr b + ab.
\]
Note also that the product $*$ can equivalently be defined in terms of the adjoint post-Lie product~\eqref{eq:rightpostlie}, as
\begin{equation}
A * B := (A_{(1)} \vartr B)A_{(2)} .\label{eq:postLie2}
\end{equation}

Recall that  $\bar{\g}:=(\mathcal{A}, \llbracket \cdot , \cdot \rrbracket)$ is related to $\g$ via \eqref{def:post-lie3} in Proposition \ref{prop:post-lie}.

\begin{Theorem}
The Hopf algebra $(\mathcal{U}(\g), *, \Delta)$ is isomorphic to the enveloping algebra $\mathcal{U}(\bar{\g})$. 
\end{Theorem}
This is essentially Theorem 3.14 in \cite{oudom2008otl}, generalized to our noncommutative setting. One first notes that the product $*$ respects the filtration of $\U(\g)$ by the length of the monomials, and that this graded product coincides with the concatenation product in each degree. The graded isomorphism $gr \mathcal{U}(\bar{\g}) \rightarrow \mathcal{U}(\g)$ results in an isomorphism $\mathcal{U}(\bar{\g}) \rightarrow (\mathcal{U}(\g), *, \Delta)$. See the proof of \cite[Theorem 3.14]{oudom2008otl}.


\subsection{Planar trees and the Grossman--Larson product}
\label{sect:GL}

We turn to the free post-Lie algebra, which as a vector space is the free Lie algebra over the set of all planar rooted trees \cite{munthe2003enumeration,munthe-kaas2013opl, vallette2007hog}. The product~\eqref{def:preLie} constructed in Section \ref{sect:UnivLieEnvAlgPostLie} corresponds to the Grossman--Larson product on planar rooted trees \cite{grossman1989has}. The dual of this product is the coproduct $\DMKW$ defined in \cite{munthe-kaas2008oth}, and we remark that the constructions of Section \ref{sect:UnivLieEnvAlgPostLie} provide a new perspective on the Hopf algebra $\HMKW$ underlying Lie--Butcher theory. 

We restate a few facts about the free post-Lie algebra. More details can be found in \cite{munthe-kaas2013opl}. The free post-Lie algebra \[\g_F = \PostLie(\ab)\] on one generator is the free Lie algebra over the set $\OT$ of planar rooted trees, equipped with a post-Lie product $\car$ called \emph{left grafting}. The first few elements of $\OT$ are 
$$
	\left\{ \begin{array}{c} 
			 \ab, \aabb,\aaabbb, \aababb, \aaaabbbb,
				\aaababbb,\aaabbabb, \aabaabbb, \aabababb,\ldots 
				\end{array}  \right\}
$$
The left grafting operation on two trees $\tau_1 \car \tau_2$ is the sum of all trees resulting from attaching the root of $\tau_1$ to all the nodes of $\tau_2$ from the left.
\begin{align}
  \AabB\,\car\aababb=\aAabBababb +\aaAabBbabb + \aabaAabBbb.
\end{align}
On brackets the operation acts as a derivation from the left, and as a difference of associators from the right. See \cite[Proposition 3.1]{munthe-kaas2013opl}. 

The post-Lie operation $\car$ has a unique extension to the universal enveloping algebra $\U(\g_F)$ by Proposition~\ref{prop:post-lieUg}. The elements of $\U(\g_F)$ are the ordered finite words in the alphabet $\OT$, including the empty word $\one$, and is called the set $\OF$ of \emph{ordered forests}.
\[
	\OF = \left\{ \begin{array}{c} 
				\one, \,\, \ab, \,\, \ab\,\ab, \,\, \aabb, \,\, \ab\,\ab\,\ab, \,\, 
				\aabb\,\ab, \,\, \ab\,\aabb, \,\, \aababb, \,\, \aaabbb, \,\, \cdots
				\end{array} \right\}.
\]
The result of left grafting a forest $\omega = \tau_1 \tau_2$ on a forest $\nu$ is a sum of words obtained by left grafting $\tau_2$ to all the nodes in $\nu$, then $\tau_1$ to all the nodes not originally from $\tau_2$. 

Using the left grafting operation one can define the important \emph{Grossman-Larson product}:
\begin{align}\label{GLprod}
\omega_1 \circ \omega_2 = \Bminus(\omega_1 \car \Bplus(\omega_2)),
\end{align}
where $B^+$ produces a tree from a forest by adding a root and $B^-$ is the inverse operation of removing a root, producing a forest.

\begin{Proposition}
In $\U(\g_F)$, the product $*$ defined in Proposition \ref{prop:postpre-lie} is the Grossman-Larson product $\circ$. 
\end{Proposition}

In \cite{munthe-kaas2008oth} the Hopf algebra $\HMKW$ underlying Lie--Butcher theory and numerical integration on Lie groups and homogeneous manifolds was introduced. It was shown that it is the dual of the Hopf algebra $(\U(\g_F), \circ, \Delta)$, and, as mentioned, the discussion in Section \ref{sect:UnivLieEnvAlgPostLie} provides an alternative point of view for $\HMKW$. 

Note that all the above constructions can be generalized to multiple generators by considering \emph{labelled} trees and forests \cite{munthe-kaas2013opl}. We refer the reader to \cite{vallette2007hog} for an operadic approach to post-Lie algebras.

Since it is established that $\circ$ and $\car$ coincide with the general post-Lie operations $*$ respectively $\tr$, we will in the sequel use the symbols $*$ and $\tr$ also on forests.


\subsection{Lie--Butcher series}
\label{ssect:LBseries}

We consider infinite series in the free post-Lie algebra, and the interpretation of these as vector fields and flows on manifolds.  Consider $\U(\g_F)$ as a vector space graded by the number of nodes in the forests, and let  
\[
	\U(\g_F)^*=\lim_{\stackrel{\longleftarrow}{n}}(\U(\g_F))_n
\]
denote the graded completion of $\U(\g_F)$, i.e.\ $\U(\g_F)^*$ consists of all infinite sums in $\U(\g_F)$ with the inverse limit topology, whose open sets are generated by sequences agreeing up to order $n$. We let the  pairing $\langle\cdot,\cdot\rangle\colon \U(\g_F)^*\times \U(\g_F)\rightarrow k$ be defined such that the ordered forests form an orthonormal set
\[
	\langle\omega,\omega'\rangle = 	
		\begin{cases}
			0 & \text{if $\omega\neq \omega'$}\\
			1 & \text{if $\omega=\omega'$.}
		\end{cases}
\]
Following~\cite{reutenauer2003free} we may identify $\U(\g_F)^*$ with the graded dual of $\U(\g_F)$. 

\begin{Definition}[Universal Lie--Butcher series] A universal LB series is an element $\alpha\in \U(\g_F)^*$. This can be written as an infinite sum 
\[
	\alpha = \sum_{\omega\in\OF}\langle\alpha,\omega\rangle\omega .
\]
\end{Definition}

Of particular importance are two special subspaces of $\U(\g_F)^*$:

\begin{Definition}\label{def:characters}
The \emph{infinitesimal characters} (primitive elements) $\ichar_F\subset \U(\g_F)^*$ and the \emph{characters} (group-like elements)
 $\Char_F\subset \U(\g_F)^*$
are defined as
\begin{align}
	\ichar_F &\coloneqq\stset{\alpha\in \U(\g_F)^*}{\Delta(\alpha) = I\otimes \alpha+\alpha\otimes I}\\
	\Char_F &\coloneqq\stset{\alpha\in \U(\g_F)^*}{\Delta(\alpha) = \alpha\otimes \alpha } ,
\end{align}
where $\Delta$ is the coshuffle coproduct.
\end{Definition}
\begin{Remark}
The name \emph{(infinitesimal) characters} originate from the identification of an infinite series $\alpha\in \U(\g_F)^*$  with  an element of the dual space, $\alpha\colon \U(\g_F)\rightarrow k, \omega\mapsto \langle \alpha,\omega\rangle$. Since the right part of the pairing is always a finte linear combination, the evaluation of the pairing is always a finite sum.  Computations on infinite series is always done by throwing the computation onto the finite right hand side, e.g.\ for $\alpha,\beta\in \U(\g_F)^*$ we compute the composition as
\[\langle\alpha*\beta,\omega\rangle = \langle \alpha\otimes \beta,\Delta_*(\omega)\rangle ,\]
where $\Delta_*\colon \U(\g_F)\rightarrow \U(\g_F)\otimes \U(\g_F)$ is the dual of the non-planar Grossman--Larson product $*$. Thus primitive elements (resp.\ group-like elements) in the completion $\U(\g_F)^*$ are infinitesimal characters (resp.\ characters) on the Hopf algebra $\HMKW$, having shuffle as product and $\Delta_*$ as coproduct. Computational formulas for $\HMKW$ can be found in~\cite{munthe-kaas2008oth}. 
\end{Remark}

Note that $\Delta$ preserves the grading and the definitions are therefore induced from the same conditions on all finite components. Note also that $\ichar_F$ is the same as the graded completion of $\g_F$. 
The concatenation product $\conc$ and the Grossman--Larson product $*$ defined on $\U(\g_F)$ extend to $\U(\g_F)^*$. We define exponentials with respect to  these products for $f\in \U(\g_F)^*$  as 
\begin{align}
	\exp^*(f) &= 1 +f + \frac12 f*f + \frac16 f*f*f +\cdots\\
 	   \exp(f) &=  1 +f + \frac12 ff  + \frac16 fff + \cdots .
\end{align}
It can be shown that these exponentials restricted to $\ichar_F$ are bijections $\exp^*$, $\exp \colon\ichar_F\rightarrow \Char_F$. The two Lie brackets $\llbracket\cdot,\cdot\rrbracket$  and $[\cdot,\cdot]$ defined on $\g_F$ extend to $\ichar_F$, turning the infinitesimal characters $\ichar_F$ into a post-Lie algebra. The characters $\Char_F$ form a group with both the products $*$ and $\conc$, with inverses $\exp^*(f)\mapsto\exp^*(-f)$ and $\exp(f)\mapsto\exp(-f)$, respectively. More explicit formulas for the inverses are  given by the antipode in $\HMKW$ and the antipode in the concatenation--deshuffle Hopf algebra~\cite{munthe-kaas2008oth}. The following characterization of the two group products will later be given concrete interpretation in terms of composition of flows on manifolds.

\begin{Proposition}The two group products $*$ and $\conc$ on $\Char_F$ are related as
\begin{equation}
	\exp(f)*\exp(g) = \exp(f)\left(\exp(f)\tr\exp(g)\right).
\end{equation}
\end{Proposition}

\begin{proof}
Since $\exp(f)$ is a character, $\Delta(\exp(f)) = \exp(f)\otimes\exp(f)$ and the result follows from~\eqref{def:preLie}.
\end{proof}

We have the following characterization of the two exponentials in terms of differential equations:

\begin{Proposition}For an infinitesimal character $f\in \ichar_F$, the curve $y(t)=\exp^*(tf)\in \Char_F$ solves the differential 
equation
\begin{equation}
\label{eq:pldiffeqn1}
	y'(t) = y(t)\left(y(t)\tr f\right), \quad y(0) = I,
\end{equation}
while the curve $z(t)=\exp(tf)\in \Char_F$ solves
\begin{equation}
\label{eq:pldiffeqn2}
	z'(t) = z(t)f, \quad z(0) = I.
\end{equation}
\end{Proposition}

\begin{proof} 
Differentiating $y(t) = \exp^*(tf)=\sum_{j=0}^\infty \frac{t^j}{j!}f^{j*}$, we get
$y'(t) = y(t)*f$. Using $\Delta(y(t)) = y(t)\otimes y(t)$ and~\eqref{def:preLie}, we obtain~\eqref{eq:pldiffeqn1}. Differentiation of $z(t)=\exp(tf)$ yields $z'(t) = z(t)f$. 
\end{proof}

Note: For $y(t)\in \Char_F$, the expression $y(t)\tr f$ is on a manifold interpreted as the parallel transport  (by the connection $\tr$) of  the vector field $f$ along the flow of $y(t)$. Thus~\eqref{eq:pldiffeqn1} is an abstract version of~\eqref{eq:lgdiffeqn}. 
The reason for the opposite ordering in the abstract setting versus the concrete manifold formulation is that the abstract group $\Char_F$ represents pullback of functions on a manifold, and hence the product $*$ maps contravariantly to composition of diffeomorphisms.
We could alternatively use the other post-Lie product $\vartr$ defined in Proposition~\ref{prop:twist}, which from~\eqref{eq:postLie2} yields $y'(t) = \left(y(t)\vartr f\right)y(t)$. However, the connection $\vartr$ on the manifold is obtained from left trivialization, in which case~\eqref{eq:lgdiffeqn} also takes the opposite form $y'(t) = y(t)\cdot f(y(t))$.

Several  characterizations of the solution operator $f\mapsto\exp^*(t f)$ are known in the literature. For most purposes it is sufficient to characterize $\exp^*(\ab)$, since the exponential of a general $tf\in \ichar_F$ can be recovered from this by  the \emph{substitution law} $\ab\mapsto tf$, see \cite{lundervold2013backward}.
The explicit
form
\begin{eqnarray*}
	\exp^*(\ab) & = &
	I + \frac{1}{2!}\ab + \frac{1}{3!}(\ab\ab+\aabb) + \frac{1}{4!}(\ab\ab\ab+\aabb\ab
	+2\ab\aabb +\aababb+\aaabbb)\\ 
	& &
	+ \frac{1}{5!}(\ab\ab\ab\ab+\aabb\ab\ab+ 2\ab\aabb\ab +3\ab\ab\aabb
	+\aababb\ab+ \aaabbb\ab+3\aabb\aabb+3\ab\aababb+3\ab\aaabbb\\ 
	& &
	+\aabababb+\aaabbabb+2\aabaabbb+\aaababbb+\aaaabbbb)
	+\frac{1}{6!}(\ab\ab\ab\ab\ab+\cdots)+\cdots
\end{eqnarray*}
is derived with recursion formulas for the coefficients in~\cite{owren1999rkm} and also discussed in~\cite{lundervold2011hao}, where the coefficients are related to \emph{non-commutative Bell polynomials}~\cite{schimming1996nbp, ebrahimi-fard2014ncb}.

The following result is of fundamental importance in many applications:

\begin{Proposition}\label{prop:exponentials}
The two exponentials $\exp^*, \exp \colon \ichar_F\rightarrow \Char_F$ are related as
\begin{equation}
	\exp^*(f) = \exp(\theta(1)) ,
\end{equation}
where $\theta(t)\in \ichar_F$ satisfies the differential equation
\begin{equation}
\label{eq:dexpinv}
	\theta'(t) = d\exp^{-1}_{\theta(t)}\left(\exp(\theta(t))\tr f\right),\quad \theta(0)=0,
\end{equation}
and where $d\exp_\theta^{-1}$ denotes the inverse left trivialized differential of the exponential map, given as an infinite series in terms of the Bernoulli numbers $B_k$ as
\begin{equation}
	d\exp^{-1}_{\theta}(\theta') = \sum_{k=0}^\infty\frac{B_k}{k!}\operatorname{ad}^{(k)}_{\theta}(\theta')
						 = \theta' + \frac12[\theta,\theta']+\frac{1}{12}[\theta,[\theta,\theta']]+ \cdots .
\end{equation}
\end{Proposition}

\begin{proof}
By differentiating $y(t) = \exp^*(tf) = \exp(\theta(t))$ and using $d/dt\exp(\theta)  = \exp(\theta)d\exp_\theta(\theta')$, we obtain
\[
	y'(t) = y(t) \left(y(t)\tr f\right) = y(t)d\exp_\theta(\theta'), 
\]
from which the result follows by inverting the linear operator $d\exp_\theta$.
\end{proof}

\begin{Remark}{\rm{The differential equation \eqref{eq:dexpinv} leads to an interesting expansion of $\theta$. A more precise description of its solution $\theta(t)$ will be given in a forthcoming article. 
First, we already stated above that
$$
	 {d\operatorname{exp}}^{-1}_{\theta(t)} := 
	 \frac{\operatorname{ad}_{\theta(t)}}{\exp(\operatorname{ad}_{\theta(t)}) -1}
	 	=\sum_{s \ge 0} \frac{B_s}{s!}\operatorname{ad}^{(s)}_{\theta(t)},
$$
where $\operatorname{ad}^{(s)}_{a}(b):=\operatorname{ad}^{(s-1)}_{a}([a,b])$. Second, $\exp\big(\theta(t)\big) \tr f = f + \sum_{i>0} \frac{1}{i!} \big(\theta(t)\big)^i\tr f$. Since $f$ is constant, we may use the ansatz $\theta(t)=\sum_{n>0} \theta_n(f)t^n$, which yields
\begin{eqnarray*}
	 {d\operatorname{exp}}^{-1}_{\theta(t)} &=& 1+ \sum_{s > 0} t^s \sum_{j=1}^{s}  \frac{B_j}{j!}  
	 \sum_{k_1 + \cdots + k_j=s \atop k_l>0} \operatorname{ad}_{\theta_{k_1}} \cdots \operatorname{ad}_{\theta_{k_j}}
\end{eqnarray*}
and
$$
	\exp\big(\theta(t)\big)\tr f = f + \sum_{i>0} \frac{t^i}{i!}  \sum_{u=1}^{i} \sum_{k_1 + \cdots + k_u=i \atop k_l>0} 
	({\theta_{k_1}} {\theta_{k_2}} \cdots {\theta_{k_u}}) \tr f .
$$
The composition of the two leads to a rather intriguing Magnus-type expansion: $\theta_1(f)=f$, $2\theta_2(f) = f \tr f$ and
\begin{eqnarray*}
	n \theta_n(f) &=& \sum_{j=1}^{n-1} \frac{1}{j!} \sum_{k_1 + \cdots +k_j=n-1 \atop k_i>0} 
	({\theta_{k_1}} {\theta_{k_2}} \cdots {\theta_{k_j}}) \tr f \\
	& & 
	+  \sum_{j=1}^{n-1} \frac{B_j}{j!} \sum_{k_1 + \cdots +k_j=n-1 \atop k_i>0} 
	{\rm ad}_{\theta_{k_1}} \cdots {\rm ad}_{\theta_{k_j}} f \nonumber\\
	& &+ \sum_{j=2}^{n-1}\bigg( \Big(  \sum_{q=1}^{j-1} \frac{B_q}{q!} \sum_{k_1 + \cdots +k_q=j-1 \atop k_i>0} 
	{\rm ad}_{\theta_{k_1}} \cdots {\rm ad}_{\theta_{k_s}}\Big) \times\\
	&& \times
	\Big( \sum_{p =1}^{n-j} \frac{1}{p!} \sum_{k_1 + \cdots +k_p=n-j \atop k_i>0} 
	({\theta_{k_1}} {\theta_{k_2}} \cdots {\theta_{k_p}}) \tr f\Big)\bigg).
\end{eqnarray*}
}}
\end{Remark}


\subsection{Vector fields and flows on manifolds} 
\label{ssect:VFflowsManifolds}

We briefly discuss post-Lie algebras appearing in differential geometry, providing concrete interpretations of the abstract notions discussed above.

A prime example of a  post-Lie algebra in differential geometry is the set of vector fields on a Lie group, with the connection given by the (right or left) Maurer--Cartan form. We will use this example to concretize the abstract operations discussed above.
 This example generalizes to the concept of \emph{post-Lie algebroids}, defined in~\cite{munthe-kaas2013opl} as  Lie algebroids equipped with a flat connection with constant torsion.

Let $(\g,[\cdot,\cdot]_\g)$ be a finite-dimensional Lie algebra with a corresponding Lie group $G$ and write $\exp_\g\colon\g\rightarrow G$ for the classical Lie exponential. For example, if $G$ is a matrix Lie group, then $\exp_\g$ denotes the matrix exponential. Vector fields on $G$ can be trivialized by either left or right multiplication. This gives rise to two adjoint post-Lie algebras. We choose right trivialization here.

Consider the vector bundle $\g\times G\rightarrow G$, the left Lie algebra action ${L}\colon \g\times G\rightarrow G$ and the \emph{anchor map} ${\lambda}\colon \g\times G\rightarrow TG$ defined as
\begin{align}
\label{eq:act1}
	{L}(V,x) &:= \exp_\g(V)\cdot x\\
	{\lambda}(V,x) &\equiv V\cdot x := \left.\frac{\partial}{\partial t}\right|_{t=0} {L}(tV,x) .\label{eq:act2}
\end{align}
Let $\g^G$ denote functions from $G$ to $\g$ and let ${\mathcal X}(G)$ denote the vector fields on $G$. A smooth $f\in \g^G$ is identified with a section of the vector bundle  $x\mapsto (f(x),x)\colon G\rightarrow \g\times G$, and the section maps to a vector field $F\in {\mathcal X}(G)$ by composition with the anchor ${\lambda}$, 
\[
	F(x) = {\lambda}(f(x),x) = f(x)\cdot x ,
\]
which we write compactly as $F = {\lambda}(f)$. In the case of Lie groups, the identification $f\mapsto F$ is  a bijection, since $f$ is recovered from $F$ by right trivialization (composition of $f$ with the right Maurer--Cartan form). For more general post-Lie algebroids over homogeneous spaces, this identification is not injective. On $\g^G$ we define a Lie bracket and a connection
\begin{align}
	[f,g](x) &:= [f(x),g(x)]_\g\\
	   f\tr g &:= {\lambda}(f)(g) \quad\Rightarrow\quad (f\tr g)(x) 
	   	      = \left.\frac{\partial}{\partial t}\right|_{t=0}g(\exp_\g(tf(x))\cdot x) .
\end{align}
For the rest of this section we put $\L:= (\g^G,[\cdot,\cdot],\tr)$.

\begin{Proposition}
$\L = (\g^G,[\cdot,\cdot],\tr)$ is  post-Lie. The anchor map ${\lambda}$ maps the bracket
$\llbracket\cdot,\cdot\rrbracket$ defined in Proposition~\ref{prop:post-lie} to the Jacobi bracket of vector fields on $G$, hence $(\g^G,\llbracket\cdot,\cdot\rrbracket,{\lambda})$ is a Lie algebroid.
\end{Proposition}

\begin{proof}
See~\cite{munthe-kaas2013opl}.
\end{proof}

Let $\U(\L)$ denote the enveloping algebra of $\L$. We can identify this with $\U(\L) = (\U(\g)^G,\conc)$, the set of maps from $G$ to $\U(\g)$ where the  product $\conc(f,g) =: fg$ is given as
\begin{equation}
	(fg)(x) = f(x) g(x)\quad \mbox{for $f,g\in\L$ and $x\in G$.}
\end{equation}
With this product $\U(\L)$ is a graded algebra and we can form the graded completion $\U(\L)^*$. 
As in Definition~\ref{def:characters} we have the infinitesimal characters $\ichar\subset \U(\L)^*$ and the characters $\Char\subset \U(\L)^*$. The infinitesimal characters represent vectorfields on $G$, while the characters represent diffeomorphisms. Note, however, that these are formal series defined in the inverse topology, and will not necessarily converge in the analytical sense. In applications, the remedy for a lack of convergence is truncation of the infinite series and estimates of exponential closeness~\cite{hairer2006gni}.

The definition of a free object implies that for any post-Lie algebra $\L$ and any assignment $\g_F\ni\ab\mapsto f\in \L$, there exists a unique post-Lie morphism $\Fel\colon \g_F\rightarrow \L$ such that $\Fel(\ab) = f$. This  extends to a unique morphism $\Fel\colon \U(\g_F)\rightarrow \U(\L)$, see~\cite{munthe-kaas2008oth,munthe-kaas2013opl}.

Therefore, for any universal Lie--Butcher series there corresponds an infinite series we call a \emph{Lie--Butcher series in $\L$}
\begin{equation}
	\Bel(\alpha)\coloneqq \sum_{\omega\in \OF} \langle\alpha,\omega\rangle \Fel(\omega) .
\end{equation}
For $\alpha\in \ichar_F$  we have $\Bel(\alpha)\in \ichar$ and for $\alpha\in \Char_F$ we have $\Bel(\alpha)\in \Char$.

The post-Lie algebra $\L$ acts  on $\F(G)$, the ring of scalar functions on $G$, via a derivation defined as
\begin{equation}
	f\tr \phi := {\lambda}(f)(\phi),\quad \mbox{for $f\in \L$, $\phi\in \F(G)$.}
\end{equation}
The derivation satisfies, for all $f,{g}\in \L$ and $\phi,\widetilde{\phi}\in \F(G)$,
\begin{align}
f\tr(\phi\widetilde{\phi}) &= (f\tr\phi)\widetilde{\phi} + \phi (f\tr\widetilde{\phi})\label{eq:lpl1}\\
[f,{g}]\tr \phi &= f\tr({g}\tr\phi)- (f\tr{g})\tr\phi - {g}\tr(f\tr\phi)+({g}\tr f)\tr\phi\label{eq:lpl2}.
\end{align}
Equation~\eqref{eq:lpl2} is equivalent to
\begin{equation}
	f\tr(g\tr\phi) - g\tr(f\tr\phi) = \llbracket f,g\rrbracket\tr\phi = (f*g-g*f)\tr\phi .
\end{equation}
The post-Lie action is extended to $\tr\colon \U(\L)\times\F(G)\rightarrow\F(G)$ as in Proposition~\ref{prop:post-lieUg}, where the rightmost elements in the equations are taken from $\F(G)$.

An important issue is the identification of the exponentials $\exp^*(f)$ and $\exp(f)$, as diffeomorphisms on the domain $G$. Since the product $*$ in $\U(\L)$ models the composition of differential operators (Lie derivations on $G$), it follows from elementary Lie theory that exponentials act as pullback on the function ring:

\begin{Proposition}\label{prop:pullback}
Let $\psi_{tf}\colon G\rightarrow G$ denote the time-$t$ flow of the vector field ${\lambda}(f)\in {\mathcal X}(G)$ and $\psi_{tf}^{\leftarrow}\phi$ the pullback of a scalar function $\phi\in \F(G)$, 
defined as $(\psi_{tf}^{\leftarrow}\phi)(x) = \phi(\psi_{tf}(x))$. Then 
\begin{equation}
\label{eq:pullback}
	\psi_{tf}^{\leftarrow}\phi = \exp^*(tf)\tr\phi .
\end{equation}
\end{Proposition}

\begin{proof}
From $f\tr\phi:= \left.\frac{\partial}{\partial t}\right|_{t=0}\psi_{tf}^{\leftarrow}\phi$ it follows that  
$\frac{\partial}{\partial t}\psi_{tf}^{\leftarrow}\phi = \psi_{tf}^{\leftarrow}(f\tr\phi)$. Iterating this, we find the Taylor series of the pullback
\begin{align}
	\psi_{tf}^{\leftarrow}\phi 	&= \phi + t f\tr \phi + \frac{t^2}{2}f\tr(f\tr\phi) + \frac{t^3}{3!}f\tr(f\tr(f\tr\phi)) + \cdots\\
						&= \phi +  t f\tr \phi + \frac{t^2}{2}(f*f)\tr\phi + \frac{t^3}{3!}(f*f*f)\tr\phi + \cdots 
						  = \exp^*(tf)\tr\phi.
\end{align}
\end{proof}

Note that the pullbacks compose in a contravariant fashion
\begin{equation}
	\psi_g^{\leftarrow}\circ\psi_f^{\leftarrow} = \left(\psi_f\circ\psi_g\right)^{\leftarrow} .
\end{equation}
This is why the symbolic differential equation~\eqref{eq:pldiffeqn2} appears in opposite order compared to the manifold equation~\eqref{eq:lgdiffeqn}.

The flow $\psi_{tf}\colon G\rightarrow G$ can be found from $\exp^*(f)\tr\phi$, by choosing $\phi$ as coordinate maps on $G$, and we obtain a well-defined action of the post-Lie characters $\exp^*(f)$ on $G$
\begin{equation}
\label{eq:charaction}
	\exp^*(tf)\cdot x := \psi_{tf}(x) \quad\mbox{for $x\in G$}.
\end{equation}
Note that this is a right action: $\exp^*(f)\cdot(\exp^*(g)\cdot x) = (\exp^*(g)*\exp^*(f))\cdot x$.

To understand the geometric difference between $\exp^*(f)$ and $\exp(f)$ we define \emph{frozen sections}. For any $f\in \g^G$ and $x_0\in G$ there exist an $f_{x_0}\in \g^G$ called $f$ \emph{frozen} at $x_0$, defined
as
\begin{equation}
\label{eq:frozen}
	f_{x_0}(x) = f(x_0) \quad\mbox{for all $x\in G$.}
\end{equation}
The freezing extends to all of $\U(\L)$ as
\begin{equation}\
	(fg)_{x_0}:= f_{x_0}g_{x_0} \quad \mbox{for $f,g\in \L$.}
\end{equation}
It is easy to verify that the freezing operation $f\mapsto f_{x_0}$ satisfies the following relations for all
$f,g\in \L\simeq\g^G$:
\begin{align}
	g\tr f_{x_0} 			&= 0\quad\mbox{for all $g\in \g^G$}	\label{eq:fr2}\\
	\left(f_{x_0}\right)_{x_1} 	&= f_{x_0}						\label{eq:fr1}\\	
	(f\tr g)_{x_0} 			&= (f_{x_0}\tr g)_{x_0}			\label{eq:fr3}\\
	[f,g]_{x_0} 			&= [f_{x_0},g_{x_0}] .			\label{eq:fr4}
\end{align}
In particular~\eqref{eq:fr2} implies that $f_{x_0}\tr f_{x_0}=0$, thus a frozen section is invariant under parallel transport by itself. Equations~\eqref{eq:fr1}-\eqref{eq:fr3} show that $f\mapsto f_{x_0}$ is projection onto the frozen sections tangent to $f$ at $x_0$, and~\eqref{eq:fr4} shows that the torsion bracket is defined fiber-wise on the sections.

Recall that a geodesic of a connection (geodesic of a covariant derivative) is a curve defined such that parallel transport of a tangent to the curve along the curve
is invariant, e.g.\ the flow of $f_{x_0}$ is a geodesic of the connection because it satisfies the infinitesimal condition of invariance under parallel transport $f_{x_0}\tr f_{x_0} = 0$.
Hence the flow of a frozen section $f_{x_0}$ is a geodesic of the connection, which is tangent to $f$ at the point $x_0$.
The concatenation exponential $\exp(f)$ models the exact flow of a geodesic tangent to $f$:

\begin{Proposition}\label{prop:plexp} 
For all $f\in \L$ we have
\begin{equation}
\label{eq:geo1}
	\exp(f)_{x_0} = \exp^*(f_{x_0})  ,
\end{equation}
the action of the character $\exp(f)$ on $G$, defined in~\eqref{eq:charaction}, is given as
\begin{equation}
\label{eq:frozenflow}
	\exp(f)\cdot x_0 = \exp_\g(f(x_0))\cdot x_0
\end{equation}
and parallel transport (pullback) of any section $A\in \U(\L)\simeq \U(\g)^G$ along the tangent geodesic to $f$ at $x_0$ is given as
\begin{equation}
	\left(\exp(f)\tr A\right)(x_0) = A\left(\exp_\g(f(x_0))\cdot x_0\right).
\end{equation}
\end{Proposition}

\begin{proof}
Since $f_{x_0}*f_{x_0} = f_{x_0}f_{x_0}+f_{x_0}\tr f_{x_0}=f_{x_0}f_{x_0}$, we obtain~\eqref{eq:geo1}. From~\eqref{eq:act1}-\eqref{eq:act2} it follows that the geodesic tangent to $f$ in $x_0$ is \[\exp(tf_{x_0})\cdot x_0=\exp_\g(tf(x_0))\cdot x_0.\] The pullback formula follows similarly to the proof of Proposition~\ref{prop:pullback}.
\end{proof}


\subsection{Vector fields and flows on $\RR^n$}
\label{ssect:VFflowsRn}

The post-Lie algebra on a Lie group described in the previous section becomes particularly simple in the euclidean case $\RR^n$. As a Lie algebra we have $G=\g=\RR^n$ and $\exp_\g(f) = f$. The canonical connection on $\RR^n$ is given as
\[
	\left(f^i\frac{\partial}{\partial x^i}\right)\tr \left(g^j\frac{\partial}{\partial x^j}\right) 
		:= f^i\frac{\partial g^j}{\partial x^i}\frac{\partial}{\partial x^j} ,
\]
which is a flat connection with zero torsion $[f,g]=0$, and hence defines a pre-Lie algebra on the set of vector fields. In this case Proposition~\ref{prop:plexp} becomes

\begin{Proposition}
The action of the character $\exp(f)$ on $\RR^n$, defined in~\eqref{eq:charaction}, is 
\begin{equation}
	\exp(f)\cdot x_0 = f(x_0)+x_0,
\end{equation}
and parallel transport of a section $A\in \U(\L)\simeq \U(\g)^G$ is given as
\begin{equation}
	\left(\exp(f)\tr A\right)(x_0) = A\left(f(x_0)+ x_0\right).
\end{equation}
\end{Proposition}


\section{Numerical integration on post-Lie algebras}
\label{sect:RK}

Lie group integrators generalize traditional numerical integrators of differential equations from vector spaces to more general manifolds.
Since the initial developments of these methods in the 1990s, it has been clear that two different Lie algebras are important for their formulation and analysis. The standard formulation of Lie group integrators is based on the full Lie algebra of all vector fields on a Lie group and the sub-algebra of (right or left) invariant vector fields. The algorithms are formulated by defining the operation of \emph{freezing} a general vector field at a given point on the group, which means replacing a general vector field by an invariant vector field tangent at the given point. Basic motions used to formulate the methods are obtained by exponentiating this `frozen' vector field. A problem with the standard formulation is that it is difficult to give an abstract algebraic meaning to the operation of freezing a general vector field at a point. It is not an operation that can be defined on an abstract post-Lie algebra, without imposing additional structure. For this reason  Lie group methods and  classical Runge--Kutta methods have so far always been formulated in a concrete setting of vector fields on a manifold (or vector space).

In the present formulation, we follow a different but equivalent approach. Recall that post-Lie algebras involve two Lie algebras, $\bar{\g}$ and $\g$, where $\g$ is not a sub-algebra of $\bar{\g}$. The two are defined over the same set, with two different Lie brackets and two different enveloping algebras. Hence we can define two different exponential maps: with respect to the product in $\U(\bar{\g})$ and with respect to the product in $\U(\g)$. We have seen that in the case of vector fields on a Lie group, exponentiation of a general vector field using the product in $\U(\g)$ is equivalent to computing the flow of a vector field frozen at a point, whereas the exponential with respect to the product in $\U(\bar{\g})$ corresponds to following the exact flow of the vector field. 

This reformulation provides a fresh and fruitful point of view. The algorithms are formulated without the process of evaluating vector fields at a given point, and without the freezing approximation. As a result the methods can be abstracted to algorithms on pre- and post-Lie algebras. The evaluation at a given point on a manifold reappears only when the abstract algorithm is interpreted in the concrete setting of flows on manifolds. The formulation and analysis of the methods, on the other hand, can be pursued entirely at an abstract algebraic level, without involving the point-evaluation operation. 
 
In the sequel we will formulate Runge--Kutta methods and Lie group integrators abstractly on pre-Lie and post-Lie algebras. As a byproduct we present a surprisingly short derivation of the classical Butcher order theory for Runge--Kutta methods on pre-Lie algebras. We hope that this will pave the way for a simplified exposition and analysis of general classes of Lie group integrators.


\subsection{Integration on post-Lie algebras}



We let $\L=(\L,[\cdot,\cdot],\tr)$ denote a post-Lie algebra, with infinitesimal characters $\ichar\subset\U(\L)^*$, characters $\Char\subset\U(\L)^*$ and exponentials $\exp,\exp^*\colon \ichar\rightarrow\Char$. We have seen in the example of vector fields over a Lie group that the map $f \mapsto \exp^*(f)$ can be interpreted as the solution operator of a differential equation with vector field $f$, whereas $\exp(f)$ is interpreted as the flow of a \emph{frozen} vector field and $\exp(f)\tr g$ is the parallel transport of $g$ along a frozen vector field. The flow and parallel transport along frozen vector fields, as well as the computation of commutators $[\cdot,\cdot]$ in $\g$, are assumed to be basic operations that can be computed exactly. The basic approximation problem of Lie group integration can be stated  in a post-Lie algebra as follows:

\begin{Problem} [The fundamental problem of numerical integration] 
Given a post-Lie algebra $\L$, approximate the exponential $\exp^*\colon \ichar\rightarrow\Char$  in terms of linear combinations in $\L$,  the bracket $[\cdot,\cdot]$ on $\L$, the exponential $\exp\colon \ichar\rightarrow\Char$, products of such exponentials, and parallel transport of the form $\exp(f)\tr g$ for $f,g\in \L$. 
\end{Problem}
 
 If $\L$ is a post-Lie algebra, a solution to this approximation problem yields a Lie group integration algorithm. If $\L$ is a pre-Lie algebra, then $[\cdot,\cdot]=0$ and solutions of the approximation problem are classical integrators.

The simplest reasonable solution to this problem is Euler's method, defined as the approximation
\begin{equation}
\label{eq:euler}
	\Psi_e(f) := \exp(f).
\end{equation}
For vector fields on a Lie group, this spells out more explicitly as
\[x_{n+1}=\Psi_e(hf)\cdot x_n = \exp(hf)\cdot x_n = \exp_\g (hf(x_n))\cdot x_n,\]
where $h\in\RR$ is the stepsize, $f\in \g^G$ a map from the Lie group $G$ to the Lie algebra $\g$, and $\exp_\g\colon \g\rightarrow G$ denotes the classical Lie exponential. For vector fields on $\RR^n$ (pre-Lie case), this is the classical Euler method
\[
	x_{n+1}=\Psi_e(hf)\cdot x_n = \exp(hf)\cdot x_n = x_n + hf(x_n), 
\]
where $f\colon \RR^n\rightarrow \RR^n$.

Since $f*f = ff+f\tr f$, we find $\exp^*(f) -\exp(f) = \frac12 f\tr f+\cdots$. Thus, with a stepsize $h \rightarrow 0$, we have $\Psi_e(hf)-\exp^*(hf) = {\mathcal O}(h^2)$, i.e.\ Euler's method is a \emph{first order approximation}. Runge--Kutta (RK) methods define higher order approximations.

\begin{Definition}[{Runge--Kutta method (RK)}] An $s$-stage RK method $\PsiRK(f)$  is defined in terms of $s^2+s$ real coefficients $\{a_{ij}\}_{i,j=1}^s$, $\{b_i\}_{i=1}^s$ as follows
\begin{align}
		 K_i &= \exp\Bigg(\sum_{j=1}^s a_{ij}K_j\Bigg)\tr f, \quad i = 1,\ldots ,s,\\
	\PsiRK(f) &= \exp\Bigg(\sum_{j=1}^s b_j K_j\Bigg),
\end{align}
where $f,K_i\in \L$ and $\PsiRK(f)\in\Char\subset {\U(\L)^*}$. If $a_{ij}$ is a strictly lower triangular matrix the method is called \emph{explicit}, otherwise $K_i$ are found by solving implicit equations.
\end{Definition}

In the pre-Lie case of vector fields on $\RR^n$ the methods $x_{n+1} = \PsiRK(hf)\cdot x_n$ become classical RK-methods:
\begin{align}
		 K_i &=hf\Bigg(x_n+\sum_{j=1}^s a_{ij}K_j\Bigg), \quad i = 1,\ldots ,s,\\
	x_{n+1} &= x_n+\sum_{j=1}^s b_j K_j,
\end{align}
Classical Runge--Kutta methods can obtain arbitrary high order, provided the coefficients $\{a_{ij}\}_{i,j=1}^s$, $\{b_i\}_{i=1}^s$ satisfy Butcher's order conditions up to the given order, see below. 

In the general post-Lie setting the methods $\PsiRK$ are in general at most second order exact, see~\cite{munthe-kaas1998rkm}. There are several ways of obtaining higher order methods in a general post-Lie algebra~\cite{munthe-kaas1998rkm,owren1999rkm,celledoni2003cfl}. As an example, the RKMK class of methods~\cite{munthe-kaas1998rkm} are of order $p$ on arbitrary post-Lie algebras, provided the coefficients satisfy the classical Butcher order conditions up to order $p$. 

\begin{Definition}[RKMK methods]\label{def:RKMK}
We define $f\mapsto \Psi_{\tiny\operatorname{RKMK}}(f)\colon \ichar\rightarrow\Char$ as:
\begin{align}
	\mbox{for\quad} i & = 1,\ldots, s \\
			     U_i &= \sum_j a_{ij}{K_j}\\
			     K_i &= d\exp^{-1}_{U_i}(\exp(U_i)\tr f)\\
	\mbox{end\quad} &\ \\
	\Psi_{\tiny\operatorname{RKMK}}(f) &:= \exp\left(\sum_i b_i{K_i}\right),
\end{align}
where $U_i, K_i\in \ichar$ and
where $d\exp^{-1}\colon \ichar\times\ichar\rightarrow \ichar$ is the inverse differential of $\exp$, given as an infinite series with Bernoulli numbers $B_k$ as coefficients
\begin{equation}
	d\exp^{-1}_{U}(V) = \sum_{k=0}^\infty\frac{B_k}{k!} \operatorname{ad}^{(k)}_{U}(V)
				  = V - \frac12[U,V]+\frac{1}{12}[U,[U,V]]+ \cdots .
\end{equation}
The expansion of $d\exp^{-1}$ can be truncated to the order of the underlying RK method.
\end{Definition}


\subsection{Pre-Lie algebras and Butcher's order theory}
\label{B-series}

In this section $\L$ is a pre-Lie algebra, and $\U(\L) = (\mathcal{S}(\g),\conc,\tr)$, where $\conc$ is now the symmetric product, written as a commutative concatenation. We briefly recall the basic definitions of $B$-series, arising from the general case of post-Lie algebras and Lie--Butcher series discussed above in the special case where $[\cdot,\cdot]\equiv 0$. We then present a short outline of Butcher's order theory for (classical) Runge--Kutta methods~\cite{butcher1963cft, butcher1972aat} in the setting of pre-Lie algebras. We remark that very recently $B$-series have been characterized geometrically as a unique and universal Taylor expansion of families of local mappings which preserve all affine symmetries between the affine spaces $\{\RR^n\}_{n=1}^\infty$, \cite{mclachlan2014b,munthe2013equivariant}.

Due to the pre-Lie relation~\eqref{def:pL}, the free pre-Lie algebra~\cite{chapoton2001pla}  is spanned by non-planar rooted trees (where the ordering of the branches is neglected). Let $T$ denote the (infinite) alphabet of non-planar rooted trees 
\[
	T = \left\{ 
		\begin{array}{c}
			\ab, \aabb,\aaabbb, \aababb, \aaaabbbb,\aaababbb,\aaabbabb, \aabababb,\ldots 
		 \end{array}\right\},
\]
and $\T=\Vect_{\RR}(T)$ the $\RR$-vector space spanned by finite linear combinations of elements in $T$. The pre-Lie product $\tr \colon \T \otimes \T \rightarrow\T$ is defined in terms of the grafting of rooted trees, e.g.
\begin{align}
  \aabb\,\tr\aababb=\aaabbababb +\aaaabbbabb + \aabaaabbbb = \aaabbababb +2\aaaabbbabb.
\end{align}
The free pre-Lie algebra (in one generator) is $\g_F = \{\T,\tr\}$.

Let $F$ denote the set of forests, i.e. all finite words in letters from $T$, where the ordering of the letters in a word is neglected. The empty word is $\one \in F$. We write $\F = \Vect_{\RR}(F)$ for the $\RR$-vector space spanned by elements in $F$. With the symmetric product of forests we have $\F = {\mathcal S}(\T) = \U(\g_F)$. The pre-Lie product $\tr$ extends uniquely from $\T$ to $\F$. Note that $F$ and $T$ are in natural 1--1 correspondence via the operations of adding and removing roots, i.e., $B^+\colon\F\rightarrow \T$ and $B^-\colon\T\rightarrow \F$.

Let $\Fc = \Lin(\F,\RR)$ be the graded completion of $\F$, consisting of all infinite $\RR$-linear combinations of trees and forests. The  pairing $\langle\cdot,\cdot\rangle\colon \Fc\times \F\rightarrow\RR$ is defined for $\omega,\omega'\in F$ such that
\begin{equation}
\label{eq:dualpair}
	\langle\omega,\omega'\rangle = \left\{\begin{array}{ll} 
								\sigma(\omega) & \mbox{if $\omega = \omega'$}\\ 0& \mbox{else},
							\end{array}\right. 
\end{equation}
where the integer $\sigma(\omega)$ counts the size of the symmetry group of the forest $\omega$. It is defined by
\begin{align}
	\sigma(\one) &= 1,\\
	\sigma(t_1t_2\cdots t_k) &=\sigma(t_1)\cdots\sigma(t_k)\mu_1!\mu_2!\cdots\mu_k!,
						\quad\mbox{for all $t_1,\ldots, t_k\in T$},\\
	     \sigma(\omega\tr \ab) &=\sigma(\omega),\quad\mbox{for all $\omega\in F$},
\end{align}
where $\mu_1,\mu_2,\ldots, \mu_k$ count the number of equal trees for each of the different shapes among $t_1,\ldots,t_k$. 
Thus the symmetry factor counts the number of planar forests which are identified by the pre-Lie relation, see~\cite{munthe-kaas2008oth} for a detailed discussion of this factor.

We define a $B$-series as an element $\alpha\in \Fc$. It can be represented as an infinite series
\begin{equation}
	\alpha = \sum_{\omega\in F} \frac{\langle \alpha, \omega\rangle}{\sigma(\omega)} \omega.
\end{equation}

In this case the infinitesimal characters and characters are defined as:

\begin{Definition}[Character] 
A $B$-series $\alpha\in \Char\subset\Fc$ is a character if and only if $\langle \alpha,\one\rangle = 1$ and for all $\omega,\omega'\in F$ 
\begin{align}
	\langle \alpha,\omega\omega' \rangle &= \langle \alpha,\omega\rangle\langle \alpha,\omega'\rangle.
\end{align}
\end{Definition}

\begin{Definition}[Infinitesimal character] 
A $B$-series $\alpha\in \ichar\subset\Fc$ is an infinitesimal character if and only if $\langle \alpha,\one\rangle = 0$ and for all $\omega,\omega'\in F\backslash\{\one\}$
\begin{align}
	\langle \alpha,\omega\omega' \rangle &=0.
\end{align}
\end{Definition}

\begin{Remark} Infinitesimal characters are expressible in terms of an infinite series in trees
\begin{equation} 
	\alpha\in \ichar \Leftrightarrow \alpha = \sum_{\tau\in T} \frac{\langle\alpha,\tau\rangle}{\sigma(\tau)}\tau ,
\end{equation} 
thus $\ichar = \Tc$, the dual of $\T$. Most authors reserve the term `$B$-series' for the infinitesimal characters. However, in that case the exponentials do not map $B$-series to $B$-series, and the distinction between the infinitesimal (vector fields) and the finite (flows and numerical methods) becomes obscured. In fact, several results in the literature actually exploit character properties of B-series, often without naming them as such, while still using the term B-series. Consequently, we choose to extend the definition of B-series to all of $\Fc$.
\end{Remark}


We will derive the classical order conditions of Runge--Kutta methods due to Butcher~\cite{butcher1963cft, butcher1972aat}. To simplify the discussion, we define operations on $s$-fold tuples. Let $\F_s:=\F\times\F\times\cdots\times\F$ and define $\ichar_s$, $\Char_s$, $\T_s$ and $\Fc_s$ similarly. Define the $s$-fold exponential $\exp_s\colon\ichar_s\rightarrow\Fc_s$, $s$-fold pairing $\langle\cdot,\cdot\rangle_s\colon \Fc_s\times\F\rightarrow\RR^s$ and $s$-fold grafting $\tr_s\colon\Fc_s\times\F\rightarrow\Fc_s$ componentwise, that is, for ${\bf K}:=(K_1,\ldots,K_s) \in \Fc_s$
\begin{align}
	\exp_s({\bf K} )&=\exp_s((K_1,\ldots,K_s)) := \left(\exp(K_1),\ldots,\exp(K_s)\right)  \\
                      \langle   {\bf K} ,\omega\rangle_s  &= \langle(K_1,\ldots,K_s),\omega\rangle_s 
                      := \left(\langle K_1,\omega\rangle,\ldots,\langle K_s,\omega\rangle\right)\\  
               {\bf K} \tr_s \omega  &= \left(K_1,\ldots,K_s\right)\tr_s \omega := (K_1\tr \omega,\ldots,K_s\tr\omega).
\end{align}
For a square matrix $A=\left(a_{ij}\right)\in \RR^{s\times s}$ and a line vector $b = \left(b_i\right)\in \RR^{1\times s}$ we define linear maps $A\colon \Fc_s\rightarrow\Fc_s$ and $b\colon \Fc_s\rightarrow {\Fc}$ as
\begin{align}
	A{\bf K}=A \left(K_1,\ldots,K_s\right) &:= \left(\sum_{j=1}^sa_{1j}K_j,\sum_{j=1}^sa_{2j}
						K_j,\ldots,\sum_{j=1}^sa_{sj}K_j\right)\\
	b{\bf K}=b\left(K_1,\ldots,K_s\right)  &:= \sum_{j=1}^s b_jK_j.
\end{align}

We now define Runge--Kutta (RK) methods using this algebraic setting. 

\begin{Definition}[Runge--Kutta method] 
A RK method applied to $\ab \in \ichar$ is 
\begin{align}
		 {\bf K} &= \exp_s(A{\bf K})\tr_s \ab\\
	\PsiRK(\ab) &= \exp\left(b{\bf K}\right),
\end{align}
where ${\bf K}\in \ichar_s$. 
\end{Definition}

This yields a map $\PsiRK\colon \ichar\rightarrow\Char$, and in particular the \emph{Runge--Kutta character} $\PsiRK(\ab)\in \Char$. We want to compare  the RK-character with the \emph{exact solution character} $\exp^*(\ab)\in \Char$, explicitly given as
\begin{align}
\label{eq:expstar}
	\exp^*(\ab) &= \sum_{k=0}^\infty\frac{1}{k!}\ab^{k*}
						= \one + \ab + \frac12 \ab * \ab +  \frac16 \ab * \ab * \ab + \cdots\\
			 & = \one+\ab+\frac12(\ab\ab+\aabb) 
			 + \frac16\left( \begin{array}{c}\ab\ab\ab+3\ab\aabb+\aababb+\aaabbb \end{array}\right)+\cdots ,
\end{align}
 
\begin{Lemma}
 The exact solution character $\exp^*(\ab)\in\Char$ is given as
\begin{equation}
\label{eq:expstar2}
	\langle\exp^*(\ab),{\tau}\rangle = \frac{1}{{\tau}!}\quad \mbox{for all ${\tau}\in T$},
\end{equation}
where ${\tau}!$ denotes the tree factorial
\begin{equation}
\label{eq:treefac}
	{\tau}! = |{\tau}| {\tau}_1!{\tau}_2!\cdots {\tau}_p! \quad\mbox{for all ${\tau} 
		  =B^+({\tau}_1{\tau}_2\cdots {\tau}_p) \in T$.}
\end{equation}
\end{Lemma}

\begin{proof} 
We prove this by induction, using $\ab*\ab = \ab\ab+\ab\tr\ab$ and the fact that $\exp^*(\ab)$ is a character. Assuming~\eqref{eq:expstar2} holds for all $|{\tau}'|<k$,  we find for  $|{\tau}|=k$, ${\tau}=B^+({\tau}_1\cdots {\tau}_p)$ that
\begin{align}
	\langle\exp^*(\ab),{\tau}\rangle &= \frac{1}{k!}\langle\ab^{k*},{\tau}\rangle 
			= \frac{1}{k!}\langle\ab^{(k-1)*}*\ab,{\tau}\rangle 
			= \frac{1}{k!}\langle\ab^{(k-1)*}\tr \ab,{\tau}\rangle\\
					   &=   \frac{1}{k} \frac{1}{(k-1)!}\langle\ab^{(k-1)*},{\tau}_1\cdots {\tau}_p\rangle 
			= \frac{1}{k}\langle\exp^*(\ab),{\tau}_1\cdots {\tau}_p\rangle\\ 
			&=  \frac{1}{k}\frac{1}{ {\tau}_1!}\cdots\frac{1}{ {\tau}_p!}. 
\end{align}
Together with $\langle\exp^*(\ab),\one\rangle=1$, the result follows.
\end{proof}

\begin{Lemma}
The Runge--Kutta character is given by the recursions
\begin{align}
	\langle\PsiRK(\ab),\tau\rangle &= b\langle {\bf K},\tau\rangle_s	\label{eq:RKrecur1}\\
	   \langle {\bf K}, \ab\rangle_s &= \een:= (1,1,\dots,1)^T\in \RR^s	\label{eq:RKrecur2}\\
	\langle {\bf K},B^+(\tau_1\tau_2\cdots\tau_p)\rangle_s &
		= A\langle {\bf K},\tau_1\rangle_s\cdot A\langle {\bf K},\tau_2\rangle_s\cdots A\langle {\bf K},\tau_p\rangle_s,\label{eq:RKrecur3}
\end{align}
where $b\colon \RR^s\rightarrow\RR$, $A\colon \RR^s\rightarrow\RR^s$ and $\ \cdot\ \colon \RR^s\times\RR^s\rightarrow \RR^s$ denotes pointwise product of vectors (known as Hadamard or Schur product).
\end{Lemma}

\begin{proof}
From the definition of the s-fold pairing we find
\begin{align}
	\langle b{\bf K},\tau\rangle &= b\langle{\bf K},\tau\rangle_s\\
	\langle A{\bf K},\tau\rangle_s&= A\langle{\bf K},\tau\rangle_s .
\end{align}
We also verify that an s-fold character ${\bf C} = (C_1,\ldots,C_s)\in \Char_s$ satisfies
\begin{equation}
	\langle {\bf C},\tau_1\tau_2\cdots\tau_p\rangle_s 
	=  \langle {\bf C},\tau_1\rangle_s\cdot \langle {\bf C},\tau_2\rangle_s\cdots \langle{\bf C},\tau_p\rangle_s,
\end{equation}
where $\ \cdot\ \colon \RR^s\times\RR^s\rightarrow \RR^s$ denotes pointwise product of vectors in $\RR^s$. For an infinitesimal character $K\in \ichar$ we have the single exponential
$\langle\exp(K),\one\rangle = 1$ and  $\langle \exp(K),\tau\rangle = \langle K,\tau\rangle$.
Thus for s-fold exponentials we have
\allowdisplaybreaks{
\begin{align}
	\langle\exp_s({\bf K}),\one\rangle_s &= \een\\
 	\langle \exp_s({\bf K}),\tau\rangle_s &= \langle K,\tau\rangle_s.
\end{align}}
Since $\langle\exp(K)\tr\ab,B^+(\tau_1\tau_2\cdots\tau_p)\rangle =  \langle\exp(K),\tau_1\tau_2\cdots\tau_p\rangle$, we find
\allowdisplaybreaks{
\begin{align}
	\langle {\bf K},B^+(\tau_1\tau_2\cdots\tau_p)\rangle_s 
		&= \langle\exp_s(A{\bf K})\tr_s\ab,B^+(\tau_1\tau_2\cdots\tau_p)\rangle_s\\
		&  = \langle\exp_s(A{\bf K}),\tau_1\tau_2\cdots\tau_p\rangle_s\\
		&= \langle\exp_s(A{\bf K}),\tau_1\rangle_s\cdots \langle\exp_s(A{\bf K}),\tau_p\rangle_s \\
		&  = A\langle {\bf K},\tau_1\rangle_s\cdots A\langle {\bf K},\tau_p\rangle_s,
\end{align}}
establishing~\eqref{eq:RKrecur3}. Equations~\eqref{eq:RKrecur1}--\eqref{eq:RKrecur2} are verified in a similar way.
\end{proof}

We have established a classical result in a pre-Lie setting:

\begin{Theorem} \label{th:rkorder} 
A Runge--Kutta method in a pre-Lie algebra has order $p$ if for all $\tau\in T$ such that $|\tau|\leq p$ we have
\begin{equation}
	\langle\PsiRK(\ab),\tau\rangle = \frac{1}{\tau!},
\end{equation}
where $\PsiRK(\ab)$ is given by ~\eqref{eq:RKrecur1}--\eqref{eq:RKrecur3}.
\end{Theorem}

The  conditions up to order 4 are (in compact and componentwise notation):
\begin{equation}
	\begin{array}{||c|c|c||} \hline
		\tau & \mbox{Compact} & \mbox{Componentwise}\\ \hline
	  	  \scalebox{0.6}{\ab} & b\een = 1& \sum_i b_i = 1\\[0.1cm]
	      \scalebox{0.6}{\aabb} & bA\een = \frac12 & \sum_{ij}b_ia_{i,j} = \frac12\\[0.1cm]
	  \scalebox{0.6}{\aababb} & b(A\een\cdot A\een)=\frac13& \sum_{ijk}b_ia_{ij}a_{ik} = \frac13\\[0.1cm]
	  \scalebox{0.6}{\aaabbb} & bA^2\een = \frac16 & \sum_{ijk}b_ia_{ij}a_{jk} = \frac16\\[0.1cm]
      \scalebox{0.6}{\aabababb} & b(A\een\cdot A\een\cdot A\een) = \frac14 & \sum_{ijkl}b_ia_{ij}a_{ik}a_{il} = \frac14\\[0.1cm]
      \scalebox{0.6}{\aabaabbb}& b(A\een\cdot A^2\een) = \frac18 & \sum_{ijkl}b_ia_{ij}a_{ik}a_{kl} = \frac18\\[0.1cm]
     \scalebox{0.6}{\aaababbb} & bA(A\een\cdot A\een) = \frac{1}{12} & \sum_{ijkl}b_ia_{ij}a_{jk}a_{jl} =\frac{1}{12} \\[0.1cm]
     \scalebox{0.6}{\aaaabbbb} & bA^3\een = \frac{1}{24} &  \sum_{ijkl}b_ia_{ij}a_{jk}a_{kl} =\frac{1}{24} \\[0.2cm]
\hline
\end{array}
\end{equation}


\subsection{Order conditions for RKMK on post-Lie algebras} 

The order conditions for Lie group methods are extensively studied in the literature~\cite{munthe-kaas1998rkm, owren1999rkm, celledoni2003cfl}. The theory can be formulated abstractly in the present post-Lie setting. However, we will not give a detailed treatment of this here, but restrict ourselves to a proof of the order conditions for the RKMK methods in Definition~\ref{def:RKMK}, following the approach of~\cite{munthe1999high} in a post-Lie setting.

\begin{Theorem}
Let $f\in \L$, where $\L$ is a post-Lie algebra. If $\{a_{ij}\}_{i,j=1}^s$ and $\{b_j\}_{j=1}^s$ satisfy the classical RK order conditions of Theorem~\ref{th:rkorder} up to order $p$, then the RKMK method in Definition~\ref{def:RKMK} is also of order $p$, satisfying
\begin{equation}
\label{eq:rkmkorder}
	\Psi_{\tiny\operatorname{RKMK}}(hf)-\exp^*(hf) = {\mathcal O}(h^{p+1}). 
\end{equation}
\end{Theorem}

\begin{proof}Let $\ichar$ and $\Char$ be the infinitesimal characters and characters of the Hopf algebra of $\L$. Given $f\in \ichar$, define a vector field $F\colon\ichar\rightarrow\ichar$ as
\begin{equation}
	F(\theta) = d\exp^{-1}_{\theta}(\exp(\theta)\tr f).
\end{equation}
According to Proposition~\ref{prop:exponentials} we have $\exp^*(f) = \exp(\theta(1))$, where $\theta(t)\in \ichar$ satisfies the differential equation $\theta(0) = 0$, $\theta'(t) = F(\theta(t))$. Note that $\ichar$ is a vector space, and is therefore also naturally a pre-Lie algebra, and we can integrate the differential equation using one step of a classical s-stage RK method of order $p$, with stepsize $h=1$, starting at $\theta(0)=0$. This yields
\begin{align}
  U_i &= \sum_{j=1}^s a_{ij}K_j\\
  K_i &=F(U_i) = d\exp^{-1}_{U_i}(\exp(U_i)\tr f), \quad i = 1,\ldots ,s,\\
  \theta(1) &\approx \sum_{j=1}^s b_j K_j,
\end{align}
which is the RKMK scheme. The verification of the order~\eqref{eq:rkmkorder} is now straightforward.
\end{proof}


\section{Acknowledgments}
The authors would like to thank the referee for helpful comments, which led to a substantial improvement of this work. We acknowledge support from~{ICMAT} and the Severo Ochoa Excelencia Program, as well as the Nils-Abel project 010-ABEL-CM-2014A, and the SPIRIT project RCN:231632 of the Norwegian Research Council. The first author is supported by a Ram\'on y Cajal research grant from the Spanish government, and acknowledges support from the spanish government under project MTM2013-46553-C3-2-P. 


\begin{thebibliography}{99}

\bibitem{burde2006left}
Burde, D.,
\newblock {\em Left-symmetric algebras, or pre-Lie algebras in geometry and physics},
\newblock Central European Journal of Mathematics, {\bf 4} (2006), 323--357.

\bibitem{butcher1963cft}
Butcher, J.C.,
\newblock {\em Coefficients for the study of Runge-Kutta integration processes},
\newblock Journal of the Australian Mathematical Society, {\bf 3} (1963), 185--201.

\bibitem{butcher1972aat}
Butcher, J.C.,
\newblock {\em An algebraic theory of integration methods},
\newblock Mathematics of Computation, {\bf 26} (1972), 79--106.

\bibitem{butcher2008nmf}
Butcher, J.C.,
\newblock {``Numerical methods for ordinary differential equations''},
\newblock John Wiley \& Sons, 2008.

\bibitem{cartier2010val}
Cartier, P.,
\newblock {\em Vinberg Algebras, Lie Groups and Combinatorics},
\newblock Quanta of Maths, {\bf 11} (2010).

\bibitem{cayley1857ott}
Cayley, A.,
\newblock {\em On the theory of the analytical forms called trees},
\newblock The London, Edinburgh, and Dublin Philosophical Magazine and Journal of Science, {\bf 13} (1857), 172--176.

\bibitem{connes1998har}
Connas, A. and D. Kreimer,
\newblock {\em Hopf Algebras, Renormalization and Noncommutative Geometry},
\newblock Commun. Math. Phys, {\bf 199} (1998), 203--242.

\bibitem{chapoton2001pla}
Chapoton, F. and M. Livernet,
\newblock {\em Pre-Lie algebras and the rooted trees operad},
\newblock International Mathematics Research Notices, {\bf 2001} (2001), 395--408.

\bibitem{celledoni2003cfl}
Celledoni, E., A. Marthinsen, and B. Owren,
\newblock {\em Commutator-free Lie group methods},
\newblock Future Generation Computer Systems, {\bf 19} (2003), 341--352.

\bibitem{ebrahimi-fard2014ncb}
Ebrahimi-Fard, K., A. Lundervold and D. Manchon,
\newblock {\em Noncommutative Bell polynomials, quasideterminants and incidence Hopf algebras},
\newblock International Journal of Algebra and Computation, {\bf 24} (2014).

\bibitem{grossman1989has}
Grossman, R. and R.G. Larson,
\newblock {\em Hopf-algebraic structure of families of trees},
\newblock Journal of Algebra, {\bf 126} (1989), 184--210.

\bibitem{hairer2006gni}
Hairer, E., C. Lubich and G. Wanner,
\newblock {``Geometric numerical integration: structure-preserving algorithms for ordinary differential equations''},
\newblock Springer Series in Computational Mathematics, {\bf 31}, 2006.


\bibitem{hoffman2003cor}
Hoffman, M.,
\newblock {\em Combinatorics of rooted trees and Hopf algebras},
\newblock Transactions of the American Mathematical Society, {\bf 355} (2003), 3795--3811.

\bibitem{iserles2000lgm}
Iserles, A., H.Z. Munthe-Kaas, S.P. N{\o}rsett and A. Zanna,
\newblock {\em Lie-group methods},
\newblock Acta Numerica 2000, {\bf 9} (2000), 215--365.

\bibitem{lundervold2011hao}
Lundervold, A. and H.Z. Munthe-Kaas,
\newblock {\em Hopf algebras of formal diffeomorphisms and numerical integration on manifolds},
\newblock Contemporary Mathematics, {\bf 539} (2011), 295--324.

\bibitem{lundervold2013backward}
Lundervold, A. and H.Z. Munthe-Kaas,
\newblock {\em Backward error analysis and the substitution law for Lie group integrators},
\newblock Foundations of Computational Mathematics, {\bf 13} (2013), 161--186.

\bibitem{manchon2011ass}
Manchon, D.,
\newblock {\em A short survey on pre-Lie algebras},
\newblock Noncommutative Geometry and Physics: Renormalisation, Motives, Index Theory, (2011), 89--102.

\bibitem{munthe-kaas1995lbt}
Munthe-Kaas, H.,
\newblock {\em Lie-Butcher theory for Runge-Kutta methods},
\newblock BIT Numerical Mathematics, {\bf 35} (1995), 572--587.

\bibitem{munthe-kaas1998rkm}
Munthe-Kaas, H.,
\newblock {\em Runge-Kutta methods on Lie groups},
\newblock BIT Numerical Mathematics, {\bf 38} (1998), 92--111.

\bibitem{munthe1999high}
Munthe-Kaas, H.,
\newblock {\em High order Runge-Kutta methods on manifolds},
\newblock Applied Numerical Mathematics, {\bf 29} (1999), 115--127.

\bibitem{munthe2003enumeration}
Munthe-Kaas, H. and Krogstad, S.,
\newblock {\em On enumeration problems in Lie--Butcher theory},
\newblock Future Generation Computer Systems, {\bf 19} (2003), 1197--1205.

\bibitem{munthe-kaas2013opl}
Munthe-Kaas, H.Z. and A. Lundervold,
\newblock {\em On post-Lie algebras, Lie--Butcher series and moving frames},
\newblock Foundations of Computational Mathematics, {\bf 13} (2013), 583--613.

\bibitem{munthe2013equivariant}
Munthe-Kaas, H.Z. and O. Verdier,
\newblock {\em Equivariant Series},
\newblock arXiv preprint arXiv:1308.5824, (2013).

\bibitem{munthe-kaas2008oth}
Munthe-Kaas, H.Z. and W.M. Wright,
\newblock {\em On the Hopf algebraic structure of Lie group integrators},
\newblock Foundations of Computational Mathematics, {\bf 8} (2008), 227--257.

\bibitem{mclachlan2014b}
McLachlan, R.I., K. Modin, H.Z. Munthe-Kaas and O. Verdier,
\newblock {\em B-series methods are exactly the local, affine equivariant methods},
\newblock arXiv preprint arXiv:1409.1019, (2014).

\bibitem{murua2006tha}
Murua, A.,
\newblock {\em The Hopf algebra of rooted trees, free Lie algebras, and Lie series},
\newblock Foundations of Computational Mathematics, {\bf 6} (2006), 387--426.

\bibitem{oudom2008otl}
Oudom, J.-M. and D. Guin,
\newblock {\em On the Lie enveloping algebra of a pre-Lie algebra},
\newblock Journal of K-theory: K-theory and its Applications to Algebra, Geometry, and Topology, {\bf 2} (2008), 147--167.

\bibitem{owren1999rkm}
Owren, B. and A. Marthinsen,
\newblock {\em Runge-Kutta methods adapted to manifolds and based on rigid frames},
\newblock BIT Numerical Mathematics, {\bf 39} (1999), 116--142.

\bibitem{reutenauer2003free}
Reutenauer, C.,
\newblock {\em Free Lie algebras},
\newblock Handbook of algebra, {\bf 3} (2004), 887--903.

\bibitem{schimming1996nbp}
Schimming, R. and S.Z. Rida,
\newblock {\em Noncommutative Bell polynomials},
\newblock International Journal of Algebra and Computation, {\bf 6} (1996), 635--644.

\bibitem{vallette2007hog}
Vallette, B.,
\newblock {\em Homology of generalized partition posets},
\newblock Journal of Pure and Applied Algebra, {\bf 208} (2007), 699--725.



\end{thebibliography}

\end{document}